\documentclass{amsproc}

\newtheorem{Theorem}{Theorem}

\newtheorem{theorem}{Theorem}[section]
\newtheorem{lemma}[theorem]{Lemma}
\newtheorem{proposition}[theorem]{Proposition}

\theoremstyle{definition}
\newtheorem{definition}[theorem]{Definition}

\theoremstyle{remark}
\newtheorem{remark}[theorem]{Remark}

\numberwithin{equation}{section}

\newcommand{\abs}[1]{\lvert#1\rvert}
\newcommand{\Abs}[1]{\Vert\lvert#1\rvert\Vert}

\newcommand{\R}{{\mathbb R}}
\newcommand{\Z}{{\mathbb Z}}

\newcommand{\N}{{\mathbb N}}

\newcommand{\Fix}{{\rm Fix}}

\begin{document}

\title[Dense properties of the circle diffeomorphisms]
{Dense properties of the space of the circle diffeomorphisms 
with a Liouville rotation number}


\author{Shigenori Matsumoto}
\address{Department of Mathematics, College of
Science and Technology, Nihon University, 1-8-14 Kanda, Surugadai,
Chiyoda-ku, Tokyo, 101-8308 Japan
}
\email{matsumo@math.cst.nihon-u.ac.jp
}
\thanks{The author is partially supported by Grant-in-Aid for
Scientific Research (C) No.\ 20540096.}


\date{\today }

\begin{abstract}
Given any Liouville number $\alpha$, it is shown that various subspaces 
are $C^\infty$-dense in
the space of the orientation preserving $C^\infty$ diffeomorphisms of
the circle with rotation number $\alpha$.
\end{abstract}

\maketitle

\section{Introduction}
Denote by $F$ the group of the orientation preserving $C^\infty$
diffeomorphisms of the cirlce. For $\alpha\in\R/\Z$, denote by $F_\alpha$
the subspace of $F$ consisting of all the diffeomorphisms with rotation
number $\alpha$. 
If $\alpha$ is irrational,
 a famous theorem of A. Denjoy \cite{D}
states that any element $f$ of $F_\alpha$ is conjugate to the
rotation by $\alpha$, denoted $R_\alpha$, by an orientation preserving
 homeomorphism. Precisely, there is a unique orientation preserving
homeomorphism $H_f$ of $S^1$ such that
$$
f=H_fR_\alpha H_f^{-1}\ \mbox{ and }\ H_f(0)=0.$$

Now denote by $O_\alpha$ the subspace of $F_\alpha$ 
of all the diffeomorphisms for which $H_f$ are $C^\infty$-diffeomorphisms.
In \cite{Y1}, J.-C. Yoccoz showed that $O_\alpha=F_\alpha$ if $\alpha$
is a non-Liouville number. He also obtained (\cite{Y2})
the following $C^\infty$-density 
result even for a Liouville number $\alpha$.

\begin{theorem} \label{t1}
For any irrational number $\alpha$, the space $O_\alpha$ 
is $C^\infty$-dense in $F_\alpha$.
\end{theorem}

Henceforth in this paper we assume $\alpha$ is a Liouville number.
Our purpose is to show that subspaces of
$F_\alpha$ which are defined by intermediate regularities of $H_f$ are
also $C^\infty$-dense in $F_\alpha$. But first of all let us state
fundamental facts about them.

\begin{remark} (\cite{H}, Chapter IV, Sect.\ 6) \label{r3}

If $H_f$
is a bi-Lipschitz homeomorphism, then it is a $C^1$-diffeomorphism.
\end{remark}

The unique $f$-invariant probability measure $\mu_f$ is given
by
$$\mu_f=(H_f)_*m,$$
where $m$ is the Lebesgue measure on $S^1$.
The uniqueness of $\mu_f$, together with the ergodicity
of $m$ w.\ r.\ t.\ $f$ (\cite{H}, p.86), implies that either $\mu_f$ is
equivalent to $m$ or singular to $m$.
From this follows easily the following remark.

\begin{remark}

(i) If $\mu_f$ is equivalent to $m$,
then the homeomorphism $H_f$ maps any Lebesgue null set to a null set.

(ii) If $\mu_f$ is singular to $m$, then
$H_f$ maps some Lebesgue null set to
a conull set.
\end{remark}

In case (i), the conjugacy $H_f$ is called {\em absolutely continuous},
and in case (ii) {\em singular}.

Besides this distinction, there are H\"older conditions.
Recall that a homeomorphism $H:S^1\to S^1$ is called $\beta$-H\"older
($\beta\in(0,1]$) if there is a constant $C>0$
such that for any $x,y\in S^1$,
$$
\abs{H(x)-H(y)}\leq C\abs{x-y}^\beta.$$

Let us define subspaces of $F_\alpha$ according to the regularity of
$H_f$ w.\ r.\ t.\ these two measures. But we shall omit the subspace
corresponding to the classes between Lipschitz and $C^1$, because it
is empty by Remark \ref{r3}.
We shall list the name of the subspaces of $F_\alpha$ and the corresponding
properties of $H_f$. 

\begin{definition} \label{d1}
Let $\beta\in(0,1)$ and $k\in\N$.

\smallskip
$\bullet$ $G_{0,{\rm sing}}$: $H_f$ is singular and
is not $d$-H\"older for any $d\in(0,1)$.

\smallskip$\bullet$ $G_{0,{\rm ac}}$: $H_f$ is absolutely continuous and 
is not $d$-H\"older for any $d\in(0,1)$.

\smallskip$\bullet$ $G_\beta$: $H_f$ is bi $\beta$-H\"older but is not
$d$-H\"older for any $d\in(\beta,1)$.

\smallskip$\bullet$ $G_{1,{\rm sing}}$: $H_f$ is singular and is bi $d$-H\"older
for any $d\in(0,1)$.

\smallskip$\bullet$ $G_{1,{\rm ac}}$: $H_f$ is absolutely continuous
 and is $d$-H\"older
for any $d\in(0,1)$, but is not bi-Lipschitz.

$\bullet$ $G_k$: $H_f$ is a $C^k$ diffeomorphism, but is not a $C^{k+1}$
diffeomorphism.
\end{definition}

The subspaces in Definition \ref{d1} are mutually disjoint.
Our main theorem is the following.

\begin{Theorem} \label{T}
For any Liouville number $\alpha$, the subspaces 
defined in Definition \ref{d1} are $C^\infty$-dense
in $F_\alpha$.
\end{Theorem}

We are not going to prove the $C^\infty$-density of $G_{0,{\rm sing}}$,
which is already done in \cite{M}. Moreover it is shown there that the subspace
of $f\in F_\alpha$ such that the Hausdorff dimension of $\mu_f$ is $0$
contains a countable union of $C^0$-open and
$C^\infty$-dense subsets of $F_\alpha$.
On the other hand it is known (\cite{Y3},\cite{S}) that
the nullity of the Hausdorff dimension of $\mu_f$ implies
that $\mu_f$ is singular and that $H^{-1}_f$ is not $d$-H\"older
for any $d\in(0,1)$. Therefore we have:

\begin{remark} (\cite{M})\label{r1}
For any Liouville number $\alpha$, the space $G_{0,{\rm sing}}$ contains a countable intersection of $C^0$-open
and $C^\infty$-dense subspaces, and in particular is residual in the
$C^r$-topology ($r=0,1,\cdots,\infty$).
\end{remark}

\begin{remark}
The subspace $G_\beta$ can be divided into $G_{\beta,{\rm sing}}$
and $G_{\beta,{\rm ac}}$, and the both subspaces seem to be
$C^\infty$-dense. But we do not pursue this problem.
\end{remark}

\smallskip

Theorem \ref{T} is proved by the method of
fast approximation by conjugacy with estimate, developed in \cite{FS}.
Moreover all the cases 
can be treated by a more or less unified fashion.
Sect.\ 2 is devoted to the exposition of this method. 
Each space in Definition \ref{d1} is dealt with separately
in further sections.

The author does not claim much originality of the arguments here. For
example, V. Sadovskaya (\cite{S}) showed that
for any Liouville number $\alpha$ and $d\in[0,1]$
there exists a diffeomorphism
$f\in F_\alpha$ for which the Hausdorff dimension of $\mu_f$
is $d$.
A 
modification of her argument 
yields a proof for the density of $G_\beta$. 
The space $G_k$ is dealt with by  elaborations
of arguments found in
\cite{KH}, 12.6. 

But Theorem \ref{T} seems to be worth recording,
while the author cannot find it in the literature. This determines
him to write down this paper.

\section{Fast approximation method}

\subsection{Overview}

We assume throughout that $0<\alpha<1$ is
a Liouville number, i.\ e.\ for any $N\in\N$ and $\epsilon>0$,
there is $p/q$ ($p,q\in\N$, $(p,q)=1$)
such that
\begin{equation}\label{e11}
\abs{\alpha-p/q}<\epsilon q^{-N}.
\end{equation}

Let $G$ be one of the subsets of diffeomorphisms of $F_\alpha$ 
defined in Definition \ref{d1}. To show that $G$
is $C^\infty$-dense in $F_\alpha$, it suffices to establish
the following proposition. We denote by $d_r$ the 
$C^r$-metric on $F$, to be defined precisely later.

\begin{proposition} \label{p1}
For any $r\in\N$, there exists $f\in G$ such that 
$d_r(f,R_\alpha)<2^{-r}$.
\end{proposition}

In fact Proposition \ref{p1} implies that $R_\alpha$ is contained
in the $C^r$-closure of $G$ for any $r$, and therefore 
in the $C^\infty$-closure of $G$,
denoted ${\rm Cl}_{\infty}(G)$. 
Since ${\rm Cl}_\infty(G)$ is invariant by the conjugation by any $h\in F$,
it follows that 
$$
O_\alpha\subset {\rm Cl}_\infty(G).$$
Then by virtue of Theorem \ref{t1}, we have
$$
F_\alpha={\rm Cl}_\infty(G).$$

\bigskip

To prove Proposition \ref{p1} for $G$, we will actually show the next
proposition.

\begin{proposition}\label{p2}
For any $r\in\N$, there are sequences $\alpha_n=p_n/q_n$ converging
to $\alpha$ and $h_n\in F$ ($n\in\N$) such that the following {\rm (i), (ii)} 
and {\rm (iii)} hold.
Define $H_0={\rm Id}$, $f_0=R_{\alpha_1}$, and for any $n\in\N$
$$
H_n=h_1\cdots h_n\ \mbox{ and }\ f_n=H_n R_{\alpha_{n+1}}H_n^{-1}.$$

\noindent
{\rm (i)} $R_{\alpha_n}$ commutes with $h_n$. 
\\
{\rm (ii)} $H_n^{-1}$ converges uniformly to a homemorphism $H^{-1}$
and its inverse $H=H_f$ satisfies the correspoing properties for $G$
listed in Definition \ref{d1}.  
\\
{\rm (iii)}
$$\abs{\alpha-\alpha_1}<2^{-r-1},\ \mbox{and}\ \
d_{n+r}(f_{n-1},f_{n})<2^{-n-r-1},\ \ \forall n\geq1.$$
\end{proposition}

Notice that (iii) implies that the limit $f$ of $f_n$
is a $C^\infty$ diffeomorphism.
It also satisfies
$d_r(f,R_\alpha)<2^{-r}$.
On the other hand the limits $H^{-1}$ and $f$ satisfy
$
H^{-1}f=R_\alpha H^{-1}$.
That is, $$f=H R_\alpha H^{-1}.
$$
This shows that Proposition \ref{p1} is satisfied by $f\in F_\alpha$.

Condition (i) is useful to establish (iii), since
then
\begin{eqnarray*}
f_{n-1}-f_n&=&H_nR_{\alpha_n}H_n^{-1}-H_nR_{\alpha_{n+1}}H_n^{-1},\\
f_{n-1}^{-1}-f_n^{-1}&=&H_nR_{-\alpha_n}H_n^{-1}-H_nR_{-\alpha_{n+1}}H_n^{-1},
\end{eqnarray*}

and these can be estimated using Lemma \ref{l2} below.

\subsection{Estimates}
Next we shall summerize inequalities needed to establish Proposition
\ref{p2}. All we need are polynomial type estimates
whose degree and coefficients can be arbitrarily large.
The inequalities below are sometimes far from being optimal.

For a $C^\infty$ function $\varphi$ on $S^1$, we define 
as usual the $C^r$ norm
$\Vert\varphi\Vert_r$ ($0\leq r<\infty$) by
$$
\Vert\varphi\Vert_r=\max_{0\leq i\leq r}\sup_{x\in S^1}\abs{\varphi^{(i)}(x)}.$$

For $f,g\in F$, define
\begin{eqnarray*}
\Abs{f}_r&=&\max\{\Vert f-{\rm id}\Vert_r,\ \Vert f^{-1}-{\rm id}\Vert_r,1\},\\
d_r(f,g)&=&\max\{\Vert f-g\Vert_r,\ \Vert f^{-1}-g^{-1}\Vert_r\}.
\end{eqnarray*}

The term $\Abs{f}_r$ is used to show that $f$ is not so big in the
$C^r$-topology. On the other hand $d_r(f,g)$ is useful for showing
$f$ and $g$ are near in the $C^r$-topology. 
The following 
inequality follows easily
from the Fa\`a di Bruno formula (\cite{H}, p.42 or \cite{S}).
We have included $1$ in the definition of $\Abs{f}_r$ 
in order to have
$$ \Abs{f}_r^i\leq\Abs{f}_r^r\ \ {\rm if}\ \ i\leq r.$$
{\em Below we denote by C(r) an arbitrary constant which depends only
on $r$.}

\begin{lemma} \label{l1}
For $f,g\in F$ we have
\begin{eqnarray*}
\Vert fg-g\Vert_r&\leq& C(r)\Vert f-{\rm Id}\Vert_r\,
\Abs{g}^r_r,\\
\Abs{fg}_r&\leq & C(r)\,\Abs{f}_r^r\,\Abs{g}_r^r.
\end{eqnarray*}

\qed
\end{lemma}

The next lemma can be found as Lemma 5.6 of \cite{FS} or
as Lemma 3.2 of \cite{S}.

\begin{lemma} \label{l2}
For $H\in F$ and $\alpha,\beta\in\R/\Z$,
$$
d_r(HR_\alpha H^{-1},HR_\beta H^{-1})\leq
C(r)\,\Abs{H}_{r+1}^{r+1}\,\abs{\alpha-\beta}.
$$
\qed
\end{lemma}

For $q\in\N$, denote by $\pi_q:S^1\to S^1$ the cyclic $q$-fold covering map.

\begin{lemma} \label{l3}
Let $h$ be a lift of $\hat h\in F$ by $\pi_q$ and assume $\Fix(h)\neq\emptyset$. 
Then we have for any $r\geq0$
\begin{eqnarray*}
\Vert h-{\rm Id}\Vert_r&=&\Vert \hat h-{\rm Id}\Vert_r\,q^{r-1}, \\
\Abs{h}_r&\leq&\Abs{\hat h}_r\,q^{r-1}.
\end{eqnarray*}
\end{lemma}

{\sc Proof}. Just notice that a lift $\tilde h$ of $h$ to $\R$ is 
the conjugate of a lift $\tilde {\hat h}$ of $\hat h$ by a homothety by $q$,
i.\ e.\ $\tilde h(x)=q^{-1}\tilde{\hat h}(qx)$.
\qed

\bigskip
Let us explain more concretely the way to construct the homeomorphisms
$h_n$ and the rationals $\alpha_n=p_n/q_n$. Essentially there are
two methods.  

\subsection{Method I}
Here we shall expose a method whose idea is easy to understand
and applicable to $G_{1,{\rm sing}}$ and $G_{1,{\rm ac}}$.
First of all we construct
beforehand a sequence of diffeomorphisms $\hat h_n\in F$ with nonempty
fixed point set; ${\rm Fix}(\hat h_n)\neq\emptyset$. Next we choose a
sequence of rationals $\alpha_n=p_n/q_n$ in a way to be explained below, and
set $h_n$ to be the lift of $\hat h_n$ by the cyclic $q_n$-fold covering
map such that ${\rm Fix}(h_n)\neq\emptyset$.
Then condition (i) of Proposition \ref{p2} is automatically satisfied. 

We always choose the rationals $\alpha_n$ so as to satisfy
$$
\abs{\alpha_{n+1}-\alpha}<\abs{\alpha_{n}-\alpha},\ \ \forall n\in\N.$$
Therfore we have
$$
\abs{\alpha_n-\alpha_{n+1}}\leq 2\abs{\alpha-\alpha_n}.$$

We shall discuss how
to define $\alpha_n$ to garantee condition (iii). 
It is by induction on $n$.
{\em Here we denote any constant which
depends on $r$, $\hat h_i$ ($1\leq i\leq n$) and 
$\alpha_1,\cdots, \alpha_{n-1}$ by $C(n,r)$. 
Thus $C(n,r)$ is any constant depending only on the innitial data about
$h_i$
and the previous step of the induction.
We also denote any positive
integer
which depends on $n$ and $r$ by $N(n,r)$}.

By Lemma \ref{l3}, we have for any $1\leq i< n$,
\begin{equation}\label{e1}
\Abs{h_i}_{n+r+1}\leq \Abs
{\hat h_i}_{n+r+1}q_i^{n+r}=C(n,r),
\end{equation}

and
\begin{equation}\label{e2}
\Abs{h_n}_{n+r+1}\leq \Abs
{\hat h_n}_{n+r+1}q_n^{n+r}=C(n,r)q_n^{N(n,r)}.
\end{equation}
Of course the two $C(n,r)$'s in (\ref{e1}) and (\ref{e2})
are different.
Now we obtain inductively using Lemma \ref{l1} that
\begin{equation}\label{e3}
\Abs{H_n}_{n+r+1}\leq C(n,r)q_n^{N(n,r)}.
\end{equation}

The terms $C(n,r)$ and $N(n,r)$ in (\ref{e3}) are computed
from (\ref{e1}) and (\ref{e2})
by applying Lemma \ref{l1} successively.
Then by Lemma \ref{l2} and (\ref{e3}),
\begin{eqnarray}\label{e4}
d_{n+r}(f_{n-1},f_n)
&=&d_{n+r}(H_nR_{\alpha_n}H_n^{-1},H_nR_{\alpha_{n+1}}H_n^{-1})
\\
&\leq&C(n,r)\,q_n^{N(n,r)}\abs{\alpha_n-\alpha_{n+1}}
\\
&\leq&C(n,r)\,q_n^{N(n,r)}\abs{\alpha-\alpha_n},\nonumber
\end{eqnarray}
for some other $C(n,r)$ and $N(n,r)$.

In order to obtain (iii) of Proposition \ref{p2}, 
the rational $\alpha_n=p/q$ have only to satisfy
$$C(n,r)q^{N(n,r)}\abs{\alpha-p/q}<2^{-n-r-1},
$$
that is,
\begin{equation}\label{e5}
\abs{\alpha-p/q}<2^{-n-r-1}C(n,r)^{-1}q^{-N(n,r)}.
\end{equation}

The terms
$$
\epsilon=2^{-n-r-1}C(n,r)^{-1}\ \mbox{ and }\
N=N(n,r)$$
are already determined beforehand or by the previous step of the induction.
Since $\alpha$ is Liouville, there exists a rational $p/q$
which satisfies (\ref{e11}) for these values of $\epsilon$
and $N$.
Setting it $p_n/q_n$, we establish (iii) for the
$n$-th step of the induction.

In fact there are infinitely many choices of $p_n/q_n$.  
This enables us to assume more.

\begin{remark}\label{r4}
The denominator $q_n$ can be chosen as large as we want, compared
with anything defined before the $n$-th step.
\end{remark}

To ensure the convergence of $H_n^{-1}$, notice that
$$
\Vert H_n^{-1}-H_{n-1}^{-1}\Vert_0=\Vert(h_{n}^{-1}-{\rm Id})\circ H_{n-1}^{-1}\Vert_0
=\Vert h_{n}^{-1}-{\rm Id}\Vert_0\leq q_{n}^{-1},$$
where the last inequality follows from the fact that $h_{n}$ is
a lift of some homeomorphism $\hat h_{n}$ by the  $q_{n}$-fold cyclic
covering map such that ${\rm Fix}(h_n)\neq\emptyset$, that is,
$h_n$ is a {\em periodic oscillation of period} $q_n^{-1}$.

Therefore if we choose $q_n$ to grow rapidly (which is possible by
Remark \ref{r4}), then
$H_n^{-1}$ converges to a continuous map $H^{-1}:S^1\to S^1$.
But since  $H^{-1}$, being a uniform limit of homeomorphisms,
is monotone and 
$H^{-1}f=R_\alpha H^{-1}$, one can show that $H^{-1}$
is in fact a homeomorphism.

\subsection{A variant of Method I}

To show the $C^\infty$ density of $G_{1,{\rm sing}}$, a slightly more 
complicated construction is convenient.
Instead of setting $h_n$ to be a lift of $\hat h_n$ by the
cyclic $q_n$-fold covering map, we can set $h_n$ to be the lift by the
cyclic $Q_n$ covering, where $Q_n=K(n)q_n$, and $K(n)$
is a positive integer determined by the previous data,
i.\ e.\ $\hat h_1\cdots\hat h_{n-1}$ and $q_1,\cdots,q_{n-1}$. 
Also in this case, (\ref{e2}) takes the form
$$
\Abs{h_n}_{n+r+1}\leq \Abs
{\hat h_n}_{n+r+1}(K(n)q_n)^{n+r}= 
C(n,r)q_n^{N(n,r)}.
$$
So the estimates of (\ref{e3}) and (\ref{e4}) is of the same form
although $C(n,r)$ is changed.
Thus we can find a solution of (\ref{e5}).

\subsection{Method II}
In this method, we define $\alpha_n=p_n/q_n$ and $h_n\in F$
at the same time in the $n$-th step of the induction.
First we consider a one parameter family $\{\hat h_t\}$ ($t\in(0,1)$)
in $F$ such that  $\Fix(\hat h_t)\neq\emptyset$. We are mostly interested for small value of $t$.
Assume for any $r\in\N$, there is a constant $C(r)>0$ and an integer
$m(r)$ such that
\begin{equation} \label{ee1}
\Abs{\hat h_t}_r\leq C(r) t^{-m(r)}.
\end{equation}
Let $\delta_n$ ($n\in\N$) be a positive number. Define
$h_n$ to be the lift of $\hat h_{q_n^{-\delta_n}}$, where
$\alpha_n=p_n/q_n$ is going to be decided.
Even in this case 
we have by (\ref{ee1})
\begin{eqnarray*}
\Abs{h_n}_{n+r+1}&\leq&\Abs{\hat h_{q_n^{-\delta_n}}}_{n+r+1}\,q_n^{n+r}\\
&\leq& C(r)q_n^{\delta_nm(r)+n+r},
\end{eqnarray*}
that is,
the estimate of $\Abs{h_n}_{n+r+1}$
takes exactly the same form as in (\ref{e2}).
Likewise that of
$\Abs{h_i}_{n+r+1}$ ($1\leq i<n$)  
is the same as in (\ref{e1}).

Therefore the inequality to decide $\alpha_n=p_n/q_n$ also
takes the form of (\ref{e5}) and the Liouville property
of $\alpha$ enables us to choose $\alpha_n$. Then the correspoding
sequence $\{f_n\}$ converges to $f\in F_\alpha$ 
in the $C^\infty$-topology such that $d_r(f,R_\alpha)<2^{-r}$
for a given $r$, and also $H_n^{-1}$ uniformly to $H^{-1}$.

In the proof of the $C^\infty$-density of $G_k$, $\delta_n=1$ for any
$n$.
For $G_\beta$, $\delta_n$ is slightly varying.  
For $G_{0,{\rm ac}}$, we start with a 2 parameter family, but this is
explained later in the section for $G_{0,{\rm ac}}$.

\section{The space $G_{1,{\rm sing}}$}

The purpose of this section is to show that the space $G_{1,{\rm sing}}$
is $C^\infty$-dense in $F_\alpha$. We follow 2.4.
Thus the sequence $\{\hat h_n\}$ is defined in the first place.
Fix once and for all integers $k_n\in\N$ such that
$$
\prod_{i=1}^\infty(1-k_i^{-1})>0.$$
For example, $k_i=(i+1)^2.$
Let
$$\hat J_n=[0,1-k_n^{-1}]\ \mbox{ and }\ \hat I_n=[0, k_n^{-1}],$$
and let $\hat h_n\in F$ be a diffeomorphism such that $\hat
h_n(\hat J_n)=\hat I_n$
and $\hat h_n$ is linear on $\hat J_n$.
Let
$$
Q_n=(k_1\cdots k_{n-1}q_1\cdots q_{n-1})\, q_n.$$
Consider the lift $h_n$ of $\hat h_n$ by the cyclic $Q_n$-fold
covering map $\pi_{Q_n}$ such that ${\rm Fix}(h_n)\neq\emptyset$.

Following Paragraph 2.4, one can choose
$\alpha_n=p_n/q_n$ inductively so that 
$f_n$ converges in the $C^\infty$ topology to $f\in F_\alpha$
such that $d_r(f,R_\alpha)<2^{-r}$ for a given $r\in\N$,
and  $H_n^{-1}=h_n^{-1}\cdots h_1^{-1}$
 uniformly to a homeomorphism $H^{-1}$.
Let
$$J_n=\pi_{Q_n}^{-1}(\hat J_n)\ \mbox{ and }\ I_n=\pi_{Q_n}^{-1}(\hat
I_n).$$
Then we have $m(J_n)=1-k_n^{-1}$ and $m(I_n)=k_n^{-1}$,
where $m$ denotes the Lebesgue measure.
Given a finite subgoup $P$ of $S^1$, we
call an interval in $S^1$ a {\em $P$-interval} if its endpoints are contained
in $P$. Since $Q_n$ is a multiple of $k_{n-1}Q_{n-1}$, we have

\medskip\noindent
(3.1) Any component of $J_{n-1}$ is a $(Q_n^{-1}\Z)/\Z$-interval.

\medskip
Let
$$C=\bigcap_{i=1}^\infty J_i\ \mbox{ and }\ C_n=\bigcap_{i=1}^nJ_i.$$
Then we have by (3.1)
$$
m(C_n)=\prod_{i=1}^n(1-k_i^{-1})\ \mbox{ and therefore }\
m(C)>0.$$
Now since $(Q_{n+1}^{-1}\Z)/\Z$ is pointwise fixed by $h_j$ ($j>n$),
we also get by (3.1) 
$$h_j(J_n)=J_n\ \mbox{ if }\ j>n.$$
This implies that for any big $j$, 
$$H_j(J_n)=H_n(J_n)=H_{n-1}h_n(J_n)=H_{n-1}(I_n).$$
That is, 
$$J_n=H_j^{-1}H_{n-1}(I_n).$$
Therefore the uniform limit $H^{-1}$ of $H_j^{-1}$ satisfies
for any $n\in\N$
$$
J_n=H^{-1}(H_{n-1}(I_n)),\ \mbox{ where }\ H_0={\rm Id}.
$$
Thus we have 
$$H(C_n)=H(\bigcap_{i=1}^nJ_i)=\bigcap_{i=1}^nH(J_i)=
\bigcap_{i=1}^nH_{i-1}(I_i).$$
Again (3.1) and the linearity of $\hat h_n$
on $\hat J_n$ enable us to compute the measure. We have
$$m(H(C_n))=\prod_{i=1}^nk_i^{-1} \to 0\ \mbox{ and therefore }\ m(H(C))=0.$$

Now the $f$-invariant measure $\mu_f=H_*m$ satisfies
$$\mu_f(H(C))=(H_*m)(H(C))=m(H^{-1}(H(C)))=m(C)>0.$$ In summary we have

\medskip\noindent
(3.2)\ \ $\mu_f(H(C))>0
\ \mbox{ and }\ m(H(C))=0$.

\medskip
Clearly (3.2) implies that 
$\mu_f$ is not equivalent to $m$, that is,
$H$ is singular.

\bigskip
In the rest we shall show that $H^{-1}$ is $d$-H\"older for any 
$d\in(0,1)$. For this purpose we do not need any further
condition on the sequence $\{\hat h_n\}$.
We just need to assume that the sequence $\{q_n\}$ grows fast
compared with $\{\hat h_n\}$.
Concretely let $M_n>1$ be a constant which satisfy
$$
M_n^{-1}\abs{x-y}\leq\abs{H_n^{-1}(x)-H_n^{-1}(y)}\leq M_n\abs{x-y},\
\ \forall x,y\in S^1.$$

We assume:

\smallskip\noindent
(3.2) For any $d\in(0,1)$, there exists $n_0$ such that
if $n\geq n_0$ then $M_{n+1}Q_{n}^{d-1}\leq 1.$

\smallskip\noindent
(3.3)  $Q_{n}^{-1}\leq 8^{-1}M_{n-1}^{-1}Q_{n-1}^{-1}$ for any
$n\geq0$.

\smallskip\noindent
(3.4)  $Q_n^{-1}\leq2^{-1}Q_{n-1}^{-1}$ for any $n\geq0$.

\smallskip
Notice that the Lipschitz constant of $h_j$ is the same as that of
$\hat h_j$, and therefore $M_{n+1}$ depends only on 
$\hat h_1,\cdots,\hat h_{n+1}$, 
which are determined in the very beginning of the argument. 
After that,  we start the inductive step
 to choose the number $\alpha_n=p_n/q_n$.
Thus (3.2) is attained if simply we choose $q_n$ big enough 
(and therefore $Q_n$ big)
compared with the Lipschitz constants of $\hat h_1,\cdots,\hat h_{n+1}$,
which is possible by
Remark \ref{r4}. Likewise (3.3) and (3.4) can be attained
by choosing $q_n$ big compared with the previous data.

Now assume $n\geq n_0$ and $Q_{n+1}^{-1}\leq\abs{x-y}\leq Q_n^{-1}$.
Then
\begin{eqnarray*}
\abs{H_{n+1}^{-1}(x)-H_{n+1}^{-1}(y)}&\leq&
M_{n+1}\abs{x-y}=M_{n+1}\abs{x-y}^{1-d}\abs{x-y}^d\\
&\leq & M_{n+1}Q_n^{d-1}\abs{x-y}^d\leq\abs{x-y}^d.
\end{eqnarray*}

On the other hand we have
$$
\abs{H_{n+1}^{-1}(x)-H_{n+1}^{-1}(y)}\geq M_{n+1}^{-1}\abs{x-y}\geq
M_{n+1}^{-1}Q_{n+1}^{-1}.
$$
This means that the distance of the points $H_{n+1}^{-1}(x)$
and $H_{n+1}^{-1}(y)$ is big enough compared with  
the periods of periodic oscillations $h_{n+2},h_{n+3},\cdots$,
and thus compared with the distance of the points
$H_{n+1}^{-1}(x)$ and $H^{-1}(x)$.

To be precise, from (3.3) and (3.4) follows
$$\sum_{i=n}^\infty Q_{i}^{-1}\leq 2Q_n^{-1}\leq4^{-1}M_{n-1}^{-1}Q_{n-1}^{-1}.$$
Therefore
\begin{eqnarray*}
\abs{H_{n+1}^{-1}(x)-H^{-1}(x)}&\leq& \sum_{i=n+2}^\infty Q_{i}^{-1}\leq
4^{-1}M_{n+1}^{-1}Q_{n+1}^{-1}\\
\leq4^{-1}M^{-1}_{n+1}\abs{x-y}
&\leq&4^{-1}\abs{H_{n+1}^{-1}(x)-H_{n+1}^{-1}(y)},
\end{eqnarray*}
and likewise for $y$.
This implies
\begin{equation*}
\abs{H^{-1}(x)-H^{-1}(y)}\leq 2\abs{H_{n+1}^{-1}(x)-H_{n+1}^{-1}(y)}
\leq2\abs{x-y}^d.
\end{equation*}

We have shown the above for any $x$ and $y$ such
that $d(x,y)\leq Q_{n_0}^{-1}$. Clearly this is enough for showing
that $H^{-1}$ is $d$-H\"older.

The proof that $H$ is $d$-H\"older for any $d\in(0,1)$ uses the same
condition (3.2) $\sim$ (3.4). The analogous argument is omitted.

\section{The space $G_{1,{\rm ac}}$}

We start with
a sequence $\{\hat h_n\}$ in $F$. This time $\hat h_n$ is
a diffeomorphism supported on $[0,2^{-n-1}]$ such that 
$\Vert\hat h_n'\Vert_0>n$. We define $h_n$ to be the lift
of $\hat h_n$ by the $q_n$-fold cyclic covering map $\pi_{q_n}$
such that
${\rm Fix}(h_n)\neq\emptyset$. As is exposed in 2.3. Method I,
we can choose $\alpha_n=p_n/q_n$
inductively such that 
$f_n$ converges in the $C^\infty$ topology to $f\in F_\alpha$
with $d_r(f,R_\alpha)<2^{-r}$ for a given $r\in\N$,
and  $H_n^{-1}=h_n^{-1}\cdots h_1^{-1}$
 uniformly to a homeomorphism $H^{-1}$.

Just as in the last part of the previous section, $H^{\pm 1}$
is $d$-H\"older for any $d\in(0,1)$ if we choose $\{q_n\}$ to grow
fast compared with $\{\hat h_n\}$.

\bigskip
Next let us show that $H$ is absolutely continuous. Let 
$$\hat K_n=[2^{-n-1},1]\ \mbox{ and }\
K_n=\pi_{q_n}^{-1}(\hat K_n).
$$
Then $h_n$ is the identity on $K_n$.
Let $X=\cap_{n=1}^\infty K_n$. Then we have
$$
m(X)\geq 1-\sum_{n=1}^\infty m(S^1\setminus K_n)
\geq 1-\sum_{n=1}^\infty 2^{-n-1}
\geq 2^{-1}.$$
Clearly $H$ is the identity on the
positive measure set $X$.
This implies that for any Bore subset $B$,
$\mu_f(B\cap X)=m(B\cap X)$ and that
$\mu_f(X)=m(X)>0$.

Then the invariant measure $\mu_f$ must be equivalent to $m$.
For, if not,
$\mu_f$ is singular to $m$, i.\ e.\
there is a Borel subset $B$ such that $m(B)=1$
and $\mu_f(B)=0$. But then $m(B\cap X)=m(X)>0$
and $\mu_f(B\cap X)=0$. A contradiction shows that
$H$ is absolutely continuous.

\bigskip
Finally let us show that $H$ is not $C^1$, which implies that
$H$ is not bi-Lipschitz by Remark \ref{r3}.
Let 
$\hat I_n$ be a closed interval such that $\hat h_n'>n$ on $\hat I_n$,
and  $I_n=\pi_{q_n}^{-1}(\hat I_n)$.
Thus  $h_n'>n$ on $I_n$.
Now if $q_n$ grows sufficiently fast, then for each $i\in\N$ and
for any component $K_i^0$
of
$K_i$, there is a component $K_{i+1}^0$ of $K_{i+1}$ such that
$K_{i+1}^0\subset K_i^0$. This shows that there is a component
$K_{n-1}^0$ of $K_{n-1}$ which is contained in $\displaystyle
\bigcap_{i=1}^{n-1}K_i$.
Clearly $H_{n-1}={\rm Id}$ on $K_{n-1}^0$. Now if $q_{n}$ is
big enough, there is a component $I_n^0$ of $I_n$ which is
contained in $K_{n-1}^0$. On $I_n^0$, we have $H_n=h_n$.

Moreover if $q_{n+1},q_{n+2},\cdots$ are big enough
compared with the length of $I_n^0$, and grow fast, then
if we put
$$Z_n=I_n^0\cap\,\bigcap_{i=n+1}^\infty K_i,$$ we have $m(Z_n)>0$.
Notice that we can assume this for any $n$.
Clearly $H=h_n$ on $Z_n$. Now 
consider the set $Z_n^*$ of the points of density
of the positive measure set $Z_n$.
The set $Z_n^*$ is perfect in the sense that
$Z_n^*$ is contained in the derived set of
$Z_n^*$. Since $Z_n^*$ is contained in the closure of $Z_n$, 
$H=h_n$ on $Z_n^*$. Moreover since $Z_n^*$ is perfect,
$H'=h_n'$
if $H$ is $C^1$.
That is,
$H'>n$ on $Z_n^*$
if we assume $H$ is a $C^1$ diffeomorphism.
Since $n$ is arbitrary, this implies that $H$ is not a $C^{1}$
diffeomorphsim.

\section{The space $G_{\beta}$}
 
Here we follow 2.5 Method II.
Consider a $C^\infty$ bump function $\psi:\R\to [0,1]$
such that $\psi((-\infty,-4^{-1}])=\{0\}$,
$\psi([4^{-1},\infty))=\{1\}$,
and $\psi$ is strictly monotone increasing on $[-4^{-1},4^{-1}]$.
For any $t\in(0,1)$, define a one parameter family $\{\hat h_t\}$
in $F$ obtained by smoothly joining the two affine functions,
one $x\mapsto t^{-1}x$ on $[0,t]$ and the other
$x\mapsto t(x-1)+1$ on $[0,1]$.
Notice that the two functions coincide at $x=t(1+t)^{-1}$,
and the latter is transfered to $x\mapsto tx$ on $[-1,0]$.
Precisely we define $\hat h_t$ as follows.

\medskip
(i) If $x\in[-4^{-1}t,4^{-1}t]$, \ \ $\hat h_t(x)=(1-\psi(t^{-1}x))tx
+\psi(t^{-1}x)t^{-1}x$.

\medskip
(ii) If $x\in [4^{-1}t,\,\, t(1+t)^{-1}-4^{-1}t]$, \ \
$\hat h_t(x)=t^{-1}x$.

\medskip
(iii) If $x\in [t(1+t)^{-1}-4^{-1}t, \,\, t(1+t)^{-1}+4^{-1}t]$,

$\hat h_t(x)=(1-\psi(t^{-1}x-(t+1)^{-1}))
t^{-1}x
+\psi(t^{-1}x-(t+1)^{-1})(t(x-1)+1)$.

\medskip
(iv) If $x\in [t(1+t)^{-1}+4^{-1}t,\,\, 1-4^{-1}t]$, \ \  $\hat h_t(x)=t(x-1)+1$.

\medskip
Since we have used affine conjugation of the same bump function,
a routine computation shows that
$\hat h_t$ satisfies (\ref{ee1}) of Paragraph 2.5.

Fix once and for all a sequence $\beta_n\downarrow \beta$, 
and set
\begin{equation}\label{eee1}
t_n=q_n^{1-\beta_n^{-1}}.
\end{equation}

Define
$h_n$ to be the lift of $\hat h_{t_n}$ by the cyclic $q_n$-covering,
where $\alpha_n=p_n/q_n$ is to be determined by (\ref{e5}) in Sect.\ 2.

Since we can choose $q_n$ big enough compared with the previous data,
we can assume the following.

\begin{equation} \label{eee2}
t_1^{-1}\cdots t_{n-1}^{-1}q_n^{-1+\beta\beta_n^{-1}}\leq 1.
\end{equation}

Notice that (\ref{eee1}) implies that

\begin{equation} \label{eee3}
q_n^{-1}=(t_nq_n^{-1})^{\beta_n},\ \mbox{ and especially }\
q_n^{-1}<(t_nq_n^{-1})^\beta.
\end{equation}

Define $$\hat I_t=[t/4, 7t/12], $$
and $I_n=\pi_{q_n}^{-1}(\hat I_{t_n})$.
Notice that 
each component of $I_n$ has length  $3^{-1}t_nq_n^{-1}$, and $h_n$
is an affine transformation with magnification $t_n^{-1}$ on $I_n$
for any small $t_n$.

Let us show first of all that $H$ is not $d$-H\"older for any
$d\in(\beta,1)$.
If one chooses $q_n$ to grow fast, then for any component $I_i^0$
of
$I_i$, there is a component $I_{i+1}^0$ of $I_{i+1}$ such that
$I_{i+1}^0\subset I_i^0$. This shows that there is a component
$I_n^0$ of $I_n$ which is contained in $\displaystyle
\bigcap_{i=1}^{n}I_i$.
On $I_n^0$, $H_n$ is an affine transformation of magnification
$t_1^{-1}\cdots t_n^{-1}$.
Let $I_n^0=[x',y']$.

Denote $H^{(n+1)}=H_n^{-1}H$. Its inverse $(H^{(n+1)})^{-1}$
is the uniform limit of 
$$
h_{n+m}^{-1}\cdots h_{n+2}^{-1}h_{n+1}^{-1}$$
as $m\to\infty$, and $h_{n+1},h_{n+2},\cdots$
are periodic oscillations of small periods $q_{n+1}^{-1},q_{n+2}^{-1},\cdots$.

Define $x,y\in S^1$ by
$$
H^{(n+1)}(x)=x'\ \mbox{ and }\ H^{(n+1)}(y)=y'.$$
It is possible to choose $\{q_n\}$ so as to satisfy
$$
\sum_{i=n+1}^\infty q_i^{-1}\leq 3^{-1}t_nq_n^{-1}$$
for any $n\in\N$.
(Decompose this into two inequalities as in the previous section,
and ressort to Remark \ref{r4}.)
Then we have
$$
\abs{x-y}\leq 3\abs{x'-y'}.
$$ 
Now it follows that
\begin{eqnarray*}
&\abs{H(x)-H(y)}=\abs{H_n(x')-H_n(y')}=t_1^{-1}\cdots t_n^{-1}\abs{x'-y'}\\
&= t_1^{-1}\cdots t_n^{-1}\cdot 3^{-1}t_nq_n^{-1}=3^{-1}t_1^{-1}
\cdots t_{n-1}^{-1}q_n^{-1},
\end{eqnarray*}
and that
$$
\abs{x-y}^{d}\leq 3^d\abs{x'-y'}^d=3^d(3^{-1}t_nq_n^{-1})^d
=t_n^dq_n^{-d}.$$

Now a computation using (\ref{eee1}) shows that if $\beta_n<d$,
\begin{eqnarray*}
\abs{H(x)-H(y)}/\abs{x-y}^{d}&\geq& 
3^{-1}t_1^{-1}\cdots t_{n-1}^{-1}q_n^{d-1}t_n^{-d}\\
&=&3^{-1}t_1^{-1}\cdots t_{n-1}^{-1}q_n^{d-1}(q_n^{1-\beta_n^{-1}})^{-d}\\
&=&3^{-1}t_1^{-1}\cdots t_{n-1}^{-1}
q_n^{\beta_n^{-1}d-1}\\
&\geq& 3^{-1}t_1^{-1}\cdots t_{n-1}^{-1}.
\end{eqnarray*}

Since the RHS of this expression can be arbitrarily large, we have shown
that $H$ is not $d$-H\"older for any $d\in(\beta,1)$.

\bigskip
We next show that $H^{-1}$ is $\beta$-H\"older. 
Assume for some $n\in\N$,
\begin{equation} \label{eee4}
t_{n+1}q_{n+1}^{-1}\leq\abs{x-y}\leq t_nq_n^{-1}.
\end{equation}

Then we have
\begin{equation}\label{eee5}
\abs{H_n^{-1}(x)-H_n^{-1}(y)}\leq\abs{x-y}^\beta.
\end{equation}

In fact, since the  Lipschitz constant of $h_i^{\pm}$ is $t_i^{-1}$
\begin{eqnarray*}
\abs{H_n^{-1}(x)-H_n^{-1}(y)}&\leq&
t_1^{-1}\cdots t_n^{-1}\abs{x-y}
\\
&=&t_1^{-1}\cdots t_n^{-1}\abs{x-y}^{1-\beta}
\abs{x-y}^\beta\\
&\leq& t_1^{-1}\cdots t_n^{-1}(t_nq_n^{-1})^{1-\beta}
\abs{x-y}^\beta
\\ &=&t_1^{-1}\cdots t_{n-1}^{-1}t_n^{-\beta}q_n^{\beta-1}\abs{x-y}^\beta\\
&
=&t_1^{-1}\cdots t_{n-1}^{-1}(q_n^{1-\beta_n^{-1}})^{-\beta}q_n^{\beta-1}
\abs{x-y}^\beta
\\ &=&t_1^{-1}\cdots t_{n-1}^{-1}q_n^{-1+\beta\beta_n^{-1}}\abs{x-y}^\beta.
\end{eqnarray*}
Thus (\ref{eee5}) follows from (\ref{eee2}).

Now the rest of the proof is divided into two cases.

\medskip
{\sc Case 1} \ \ $\abs{H_n^{-1}(x)-H^{-1}_n(y)}\geq q_{n+1}^{-1}$.

Since $h_{n+1}^{-1}$ is a periodic osillation of period $q_{n+1}^{-1}$,
it follows from the assumption of Case 1 that
$$
\abs{H_{n+1}^{-1}(x)-H_{n+1}^{-1}(y)}\leq2\abs{H_n^{-1}(x)-H^{-1}_n(y)}.
$$
If we choose $q_{n+2},q_{n+3},\cdots$ to grow fast, then we have
$$\abs{H^{-1}(x)-H^{-1}(y)}\leq3\abs{H_n^{-1}(x)-H^{-1}_n(y)}
\leq 3\abs{x-y}^\beta,
$$
finishing the proof in this case.

\medskip
{\sc Case 2} \ \ $\abs{H_n^{-1}(x)-H^{-1}_n(y)}\leq q_{n+1}^{-1}$.

Since $h_{n+1}^{-1}$ is a periodic oscillation of period $q_{n+1}^{-1}$,
we have in this case that
$$
\abs{H_{n+1}^{-1}(x)-H^{-1}_{n+1}(y)}\leq q_{n+1}^{-1},
$$
and therefore we get using (\ref{eee3})
$$
\abs{H^{-1}(x)-H^{-1}(y)}
\leq 2q_{n+1}^{-1}\leq 2(t_{n+1}q_{n+1}^{-1})^\beta\leq2\abs{x-y}^\beta,
$$
and the proof is complete.

\bigskip
Finally let us show that $H$ is $\beta$-H\"older.
Again we assume for some $n\in\N$,
$$
t_{n+1}q_{n+1}^{-1}\leq\abs{x-y}\leq t_nq_n^{-1}.
$$
Let 
$$
x'=H^{(n+2)}(x),\ \mbox{ and }\ y'=H^{(n+2)}(y).$$
If $q_{n+2}^{-1}, q_{n+3}^{-1},\cdots$ are small enough, then we have
$$
\abs{x'-y'}\leq 2\abs{x-y}.$$

\medskip
{\sc Case 1} $\abs{x'-y'}\geq q_{n+1}^{-1}$.

If we put
\begin{equation} \label{eee6}
x''=h_{n+1}x'\ \mbox{ and }\ y''=h_{n+1}x',
\end{equation}
then by the assumption of Case 1, we have
$$
\abs{x''-y''}\leq 2\abs{x'-y'}.$$
We can show by the same computation as for $H^{-1}$ which
stems from the condition 
$\abs{x-y}\leq t_nq_n^{-1}$, that
$$
\abs{H_n(x'')-H_n(y'')}\leq 4t_1^{-1}\cdots t_{n}^{-1}\abs{x-y}
\leq 4\abs{x-y}^\beta,
$$
completing the proof in this case.

\medskip
{\sc Case 2} $\abs{x'-y'}\leq q_{n+1}^{-1}$.

As before we have $\abs{x''-y''}\leq q_{n+1}^{-1}$, where $x''$ and
$y''$ are defined by (\ref{eee6}). Now we have
\begin{eqnarray*}
&\abs{H(x)-H(y)}=\abs{H_n(x'')-H_n(y'')}\\
&\leq t_1^{-1}\cdots t_n^{-1}(q_{n+1}^{-1})^{1-\beta\beta_{n+1}^{-1}}
\abs{x''-y''}^{\beta\beta_{n+1}^{-1}}.
\end{eqnarray*}
Therefore by (\ref{eee2}), we have
$$
\abs{H(x)-H(y)}\leq\abs{x''-y''}^{\beta\beta_{n+1}^{-1}}.$$
On the other hand by (\ref{eee3}),
$$
\abs{x''-y''}\leq q_{n+1}^{-1}= (t_{n+1}q_{n+1}^{-1})^{\beta_{n+1}}
\leq \abs{x-y}^{\beta_{n+1}}.$$

The proof that $G_\beta$ is $C^\infty$-dense is now complete.

\section{The space $G_{0,{\rm ac}}$}

Given $0<s<t<1$, consider the three affine transformations:

$\bullet$ $x\mapsto s^{-1}tx$ for $x\in[0,s]$.

$\bullet$ $x\mapsto st^{-1}(x-t)+t$ for $x\in[0,t]$.

$\bullet$ $x\mapsto x$ for $x\in[t,1]$.

Define a family $\hat h_{s,t}$ in $F$ which smoothly joins the above three
transformations by  bump functions.
On the interval $\hat I_{s,t}=[s/4,7s/12]$, $\hat h_{s,t}$
is to be an affine transformation of magnification $s^{-1}t$.
We also have $\hat h_{s,t}={\rm Id}$ on the interval $\hat K_{s,t}
=[2t,1]$.
If we choose the bump functions
to be affine conjugates of the same function, as we did in Sect.\ 5, we have
\begin{equation}\label{eeee1}
\Abs{\hat h_{s,t}}_r\leq C(r)s^{-m(r)}t^{-m(r)}
\end{equation}
for some $C(r)$ and $m(r)$.

Let 
\begin{equation}\label{eeee0}
t=3^{-n-2}
\end{equation}
and $q_n\in\N$ be arbitrary. Let $h_n$ be the lift of $\hat h_{s,t}$
by the $q_n$-fold cyclic covering $\pi_{q_n}$ such that
$\Fix(h_n)\neq\emptyset$, and let $[x,y]$ be a component of
$I_n=\pi_{q_n}^{-1}(\hat I_{s,t})$. 
Also let $K_n=\pi_{q_n}^{-1}(\hat K_{s,t})$.

Then for any $d\in(0,1)$, we have
\begin{equation}\label{eeee2}
\frac{h_n(x)-h_n(y)}{(x-y)^d}=\frac{3^{-n-3}q_n^{-1}}{(3^{-1}sq_n^{-1})^d}.
\end{equation}

Setting the above ratio to be equal to $n$ and $d=n^{-1}$, we get
\begin{equation}\label{eeee3}
s=3^{-n^2-3n+1}n^{-n}q_n^{-n+1}.
\end{equation}

If we fix the value of $t$ and $s$ by (\ref{eeee0}) and (\ref{eeee3}),
then (\ref{eeee1}) implies that for any $r\in\N$
$$
\Abs{h_n}_{n+r+1}\leq C(n,r)q_n^{N(n,r)},$$
for some constants $C(n,r)$ and $N(n,r)$.
Therefore we can follow 2.5.\ Method II, and can choose
$\alpha_n=p_n/q_n$
inductively so that 
$f_n$ converges in the $C^\infty$ topology to $f\in F_\alpha$
in such a way that $d_r(f,R_\alpha)<2^{-r}$ for a given $r\in\N$,
and  $H_n^{-1}=h_n^{-1}\cdots h_1^{-1}$
 uniformly to a homeomorphism $H^{-1}$.

Just as in Sect.\ 4, if we choose $q_n$ to grow fast, then 
there is a component $I_n^0$ of $I_n$ such that
$
m(Z_n)>0
$, where
$$
Z_n=I_n^0\cap\,\bigcap_{i\neq n}K_i.
$$

For any distinct points $x,y\in Z_n$, we have
$$
\frac{H(x)-H(y)}{(x-y)^{n^{-1}}}=\frac{h_n(x)-h_n(y)}{(x-y)^{n^{-1}}}=n,$$
since we have set the value of (\ref{eeee2}) to be equal to $n$.
This shows that $H$ is not $d$-H\"older for any $d\in(0,1)$.

\bigskip
Finally the same argument as in Sect.\ 4 shows
that $H$ is absolutely continuous.

\section{The space $G_k$}

In this section we shall show that $G_k$ in Definition \ref{d1}
is $C^\infty$-dense in $F_\alpha$.
We consider the one parameter family $\{\hat h_t\}_{0<t<1}$ given
by Lemma \ref{last} below, and 
follows 2.5.\ Method II for the choice of
$\delta_n=1$ ($\forall n\in\N$), that is, we
define $h_n$ to be the lift of
$\hat h_{q_n^{-1}}$ by the cyclic $q_n$-fold covering such
that ${\rm Fix}(h_n)\neq\emptyset$.

\begin{lemma}\label{last}
Given $k\in\N$, there exists a one parameter family $\{\hat h_t\}_{0<t<1}$ in $F$ which
satisfies the following properties.

\smallskip\noindent
{\rm (i)} $\hat h_t ={\rm Id}$ in an interval $\hat K_t$ of length $\geq 
1-t^k$.

\smallskip\noindent
{\rm(ii)} $\Vert\hat h_t^{\pm}- {\rm Id}\Vert_k\leq C t^{k}$ for some constant $C>0$.

\smallskip\noindent
{\rm(iii)} $\abs{(\hat h_t^{-1})^{(k+1)}}\geq C^{-1}$ on an interval
$\hat I_t$ for some constant
$C>0$.

\smallskip\noindent
{\rm(iv)} For any $r\in\N$, $\Abs{\hat h_t}_r\leq C_r t^{-m(k,r)}$
for some constant $C_r>0$ and an integer $m(k,r)$.
\end{lemma}

The proof of this lemma is postponed until the end of this section.

By following 2.5.\ Method II, thanks to (iv) above,
we can choose $\alpha_n=p_n/q_n$
inductively so that 
$f_n$ converges in the $C^\infty$ topology to $f\in F_\alpha$
in such a way that $d_r(f,R_\alpha)<2^{-r}$ for a given $r\in\N$,
and  $H_n^{-1}=h_n^{-1}\cdots h_1^{-1}$
 uniformly to a homeomorphism $H^{-1}$.

Let us show that a bit more careful choice of
$\alpha_n$ ensures that $H^{-1}$ is a $C^k$ diffeomorphism.
Now (ii) above and Lemma \ref{l3} implies that
$$
\Vert h_n^{-1}-{\rm Id}\Vert_k=\Vert\hat h_{q_n^{-1}}^{-1
}-{\rm Id}\Vert_k\cdot q_n^{k-1}\leq C q_n^{-1}.$$
By Lemma \ref{l1} we have
$$
d_k(H_n^{-1},H_{n-1}^{-1})\leq C\Vert h_n^{-1}-{\rm Id}\Vert_k
\Abs{H_{n-1}}_k^{k}\leq Cq_n^{-1}\Abs{H_{n-1}}_k^{k}.$$
Thus if we choose $q_n$ big enough compared with $\Abs{H_{n-1}}_k^{k}$,
then we have
$$
d_k(H_n^{-1},H_{n-1}^{-1})\leq 2^{-n}.$$
This shows that the limit $H^{-1}$ is a $C^k$ diffeomorphism.

\bigskip
Let us show that $H^{-1}$ is not a $C^{k+1}$ diffeomorphism.
Let 
$$K_n=\pi_{q_n}^{-1}(\hat K_{q_n^{-1}})\
\mbox{ and }\ I_n=\pi_{q_n}^{-1}(\hat I_{q_n^{-1}}).$$
Thus $h_n={\rm Id}$ on $K_n$ and
 $\abs{(h_n^{-1})^{(k+1)}}>C^{-1}q_n^k$ on $I_n$,
the latter by (iii) above.
The rest of the proof is nearly the same as in Sect.\ 4,
where we have shown that $H$ is not $C^1$.

\bigskip
Now let us prove Lemma \ref{last}.
We begin with the following lemma found in \cite{H}, p.\ 154. 
The proof is included for the convenience of the reader.

\begin{lemma}
Let $\{\phi_a\}$ be a $C^\infty$ flow on $S^1$. 
For any $r\geq1$, there is a constant $C_r>0$ such that for any
$a\in[-1,1]$, we have
$$
C_r^{-1}\abs{a}\leq\Vert (\phi_a)^{(r)}\Vert_0\leq C_r \abs{a}.$$
\end{lemma}

{\sc Proof}. Let 
$$C'_r=\sup_{a\in[0,1]}\Vert\frac{\partial}{\partial a}\phi_a^{(r)}\Vert_0.
$$
Then since 
$\phi_0^{(r)}=0$, we have by the mean value theorem
$$\abs{a}^{-1}\Vert \phi_a^{(r)}\Vert_0\leq C'_r.
$$

To show the converse inequality, consider 
$\displaystyle X=\frac{\partial}{\partial
a}\phi_a\vert_{a=0}$,
the infinitesimal generator of the flow $\{\phi_a\}$. Then
$\displaystyle X^{(r)}=\frac{\partial}{\partial a}\phi_a^{(r)}\vert_{a=0}$
is not constantly equal to zero. 
Assume for some $x\in S^1$
$$ \abs{\frac{\partial}{\partial a}\phi_a^{(r)}(x)\vert_{a=0}}=
2C''_r>0.$$
Then, by the definition of $\partial/\partial a$, 
there is $\alpha>0$ such that if $\abs{a}\leq\alpha$,
$$
\abs{a}^{-1}\abs{\phi_a^{(r)}(x)}\geq C''_r.$$
That is,
$$\abs{a}^{-1}\Vert \phi_a^{(r)}\Vert_0\geq C''_r.$$
For $\alpha\leq\abs{a}\leq 1$, 
$\abs{a}^{-1}\Vert \phi_a^{(r)}\Vert_0$
is positive and continuous in $a$.
The proof is complete by relaxing the constant $C'_r$
and $(C''_r)^{-1}$.
\qed

\bigskip
Choose a flow $\{\phi_a\}$ whose support is contained in 
$(0,1)\subset S^1$.  For any $k\in\N$ and $t\in(0,1)$ let 
$$
\hat h_t(x)=t^k\phi_{t^{k^2}}(t^{-k}x).$$
Then it is routine to check the properties in Lemma \ref{last}.

\end{document}